\documentclass{amsart}
\usepackage{amssymb}
\usepackage{amsmath}
\usepackage{amsfonts}
\theoremstyle{plain}
\newtheorem{thm}{Theorem}[section]
\newtheorem{prop}[thm]{Proposition}

\newtheorem{cor}[thm]{Corollary}
\newtheorem{dfn}[thm]{Definition}
\newtheorem{rmk}[thm]{Remark}

\newtheorem{mthm}[thm]{Main Theorem}

\begin{document}

\begin{center}
Version of  February 4, 2008
\end{center}

\bigskip

\title
{Examples of cohomology
manifolds which are not homologically locally connected }

\author{Umed H. Karimov}
\address{Institute of Mathematics,
Academy of Sciences of Tajikistan,
Ul. Ainy $299^A$, Dushanbe 734063, Tajikistan}
\email{umed-karimov@mail.ru}

\author{Du\v san Repov\v s }
\address{Institute of Mathematics, Physics and Mechanics,
and Faculty of Education, University of Ljubljana, P.O.Box 2964,
Ljubljana 1001, Slovenia} \email{dusan.repovs@guest.arnes.si}



\begin{abstract}
Bredon has constructed a 2-dimensional compact cohomology manifold
which is not homologically locally connected, with respect to the singular homology.
In the present paper we construct infinitely many such examples
(which are in addition metrizable spaces) in all remaining dimensions
$n \ge 3$.
\end{abstract}

\maketitle \pagestyle{myheadings} \markright{ON NON-HLC COHOMOLOGY MANIFOLDS}

\section{Introduction}

We begin by fixing terminology, notations and formulating
some
elementary facts. We shall use singular homology $H_n$, singular
cohomology $H^n$ and \v Cech cohomology $\Check H^n$ groups with
integer coefficients $\mathbb{Z}$. By $HLC$ and $clc$ spaces we shall denote
 homology and cohomology locally connected spaces with respect
to singular homology and \v Cech cohomology, respectively.
The general references for  terms undefined in this paper
will be
\cite{Br}--\cite{H}, \cite{Ke}, \cite{Sp}.

\begin{dfn} {\rm(cf. \cite[p.377, Corollary 16.19]{Br}).}
A finite-dimensional cohomology locally connected space $X$ is
an $n$-dimensional cohomology manifolds ($n$-cm) if
\[
\check{H}^p(X, X \setminus \{x\}) = \left\{
\begin{array}{ll}
\mathbb{Z},  \quad &\mbox{for } p = n, \\
0,  \quad &\mbox{for } p \neq n
\end{array}\right.
\]
for all $x \in X.$
\end{dfn}

We denote the segment $[0, 1]$ by $\mathbb{I}$, the
$n-$dimensional cube by $\mathbb{I}^n,$ and
the
quotient space of a space
$X$ by its subset $B$ by
$X/B$. 
The double $d(X,B)$ of $X$ with
respect to $B$ is the quotient space of the product $X\times \{0,
1\}$, where $\{0, 1\}$ is the two point set, by identification of
the points $(x,0)$ with $(x,1)$ for every $x\in B.$

\begin{dfn} {\rm(cf. \cite[p.385, Theorem 16.27]{Br}).}
A space $X$ is said to be a cohomology manifold with boundary $B$
if $B$ is a closed nowhere dense subspace of $X$ and the double
$d(X,B)$ is a cohomology $n$-manifold.
\end{dfn}

Bredon constructed a 2-dimensional compact cohomology manifold
which is not homologically locally connected (non-$HLC$) space
\cite[p.131, Example 17.13]{Br}. Every metrizable 2-dimensional
locally compact cohomology manifold is a topological 2-manifold
and therefore it
is a $HLC$ space (see e.g. \cite[p.388, Theorem
16.32]{Br}).

The main goal
of the present paper is to construct
for all
remaining dimensions $n \geq 3,$ infinitely many $n$-dimensional
compact metrizable cohomology manifolds which are not
homologically locally connected.

\begin{mthm}
Let $M^n$ be an n-dimensional compact manifold (possibly with
boundary), $n \geq 3.$ Then there exists a 
compact
subset
$C$ 
of
the
interior int$(M^n)$ of $M^n$ such that the quotient space $M^n/C$
is a non-homologically locally connected n-dimensional metrizable
cohomology manifold.
\end{mthm}

\section{Preliminaries}

Let $G$ be any multiplicative  (in general nonabelian)
group. By
the {\sl commutator}
 of two elements $a,b\in G$ we mean the
following
product $[a,b]=a^{-1}b^{-1}ab \in  G$.
Let $G_n$ be the {\sl lower central series} of $G$ which is defined
inductively: $G_1 = G, G_{n+1} = [G_n, G],$ where $[G_n, G]$ is
the normal subgroup of $G$ generated by the following set of
commutators: $\{ [a,b]: a\in G_n, b\in G\}$.

If $F$ is a free group, the factor group $F/F_n$ is called a 
{\sl free nilpotent group}.
Let $A*B$ be the free product of groups $A$ and $B.$ 
Every
non-neutral element $x$ of $A*B$ is then
uniquely expressible in the
{\sl reduced form} as $x = u_1u_2 \cdots u_n,$ where all $u_i\in A\cup
B$ are non-neutral elements of the groups $A$ and $B$, and if $u_i
\in A$ then $u_{i+1}\in B$ (if $u_i\in B$ then $u_{i+1}\in A$,
respectively), for $i \in \{ 1,2,\dots, n-1\}$.

Following Rhemtulla \cite[p.578]{Rh}, we define for an element
$b\in B$ of order $> 2$ the mapping $\it{w}_b: A*B \to 
\mathbb{Z}$ as follows: Let $(x,b)$ denote
the multiplicity of $b$ in $x$,
i.e. the number of occurrences
of $b$ in the reduced form of $x.$ Similarly, let $(x,b^{-1})$ be
the multiplicity of $b^{-1}$ in $x.$ Write $\it{w}_b(x) = (x, b) -
(x, b^{-1}).$ For example, let $a$ be an element of group $A$
order of which is $>2$. Then
\begin{equation}\label{equation:ab}
\it{w}_b([ab,ba]) = \it{w}(b^{-1}a^{-1}a^{-1}b^{-1}ab^2a) = -2.
\end{equation}

By the {\it commutator length}  of $g\in G$, denoted by  $cl(g)$, we denote the
minimal number of the commutators of the group $G$ whose product
is equal to $g$. If such a number does not exist then we set
$cl(g) = \infty$. We also set  $cl(e) = 0$ for the neutral element
$e$ of the group $g\in G$. Obviously,
\begin{equation}\label{equation:infty}
cl(g) = \infty \ \ \mbox{if and only if} \ \ g\notin G_2,
\end{equation}
$$cl(g_1g_2) \leq cl(g_1) + cl(g_2)$$
or, equivalently (since $cl(g^{-1}) = cl(g)$)
\begin{equation}\label{equation:cl<}
cl(g_2) - cl(g_1) - cl(g_3) \leq cl(g_1g_2g_3).
\end{equation}

If $\varphi :G \to G'$ is a homomorphism of groups then for every
$g\in G,$
\begin{equation}\label{equation:cl(phi)}
cl(\varphi (g)) \leq cl(g).
\end{equation}

For any
path connected space X, the fundamental group $\pi_1(X)$
does not depend on the choice of the base point and $H_1(X)$ is
isomorphic to the factor group $\pi_1(X)/[\pi_1(X), \pi_1(X)].$
If $g$ is a loop in the space $X$ then by $[g]$ we denote
the
corresponding element of the fundamental
group $\pi_1(X)$.

For the construction of examples and proofs of their asserted
properties we shall
need some general facts:

\begin{prop}\label{prop:Series} {\rm(\cite[Proposition 3.4]{LS}).}
The lower central series $\{F_n\}$ of any
free group $F$ has a
trivial intersection, i.e. $\cap_{n=1}^{\infty}F_n = \{e\}.$
\end{prop}

\begin{prop}\label{prop:Torsions} {\rm(\cite[Theorem 1.5]{M}).}
Every free nilpotent group is torsion-free.
\end{prop}

\begin{prop}\label{prop:Generators}
{\rm(\cite[Exercise 2.4.13]{MKS}).}
Let $F$ be any
free group on
two generators $a, b.$ Then the
subgroup $F_2$ is freely generated by the commutators $[a^n,
b^m],$ where $n$ and $m$ are nonzero integers.
\end{prop}

Suppose that $b\in B$ and its order is $> 2$. The main
properties of the function $\it{w}_b$ of Rhemtulla which are
necessary for our
proofs are formulated in the following three
propositions:

\begin{prop}\label{prop:For 3} {\rm(\cite[page 579]{Rh}).}
For any elements $x\in A$ and $y\in B$
$$|\it{w}_b(xy)| \leq |\it{w}_b(x)| + |\it{w}_b(y)| + 3.$$
\end{prop}

\begin{prop}\label{prop:For 9} {\rm(\cite[page 579]{Rh}).}
For any commutator $[x, y], \ x\in A$ and $y\in B$,
$|\it{w}_b([x,y])| \leq 9.$
\end{prop}

By induction, using Propositions~\ref{prop:For 3} and
~\ref{prop:For 9}, one gets the following:

\begin{prop}\label{prop:k commutators} {\rm(\cite[page 579]{Rh}).}
For any $k$ commutators $[x_i, y_i],$ \ $x_i\in A$ and $y_i\in B$,
$\ i = 1,2,\dots k$,
$$|\it{w}_b(\prod_{i=1}^k [x_i,y_i])| \leq 12k - 3.$$
\end{prop}

\begin{prop}\label{prop:Sup}
Let $A*B$ be the free product of groups $A$ and $B$. Let $a\in A$
and $b \in B$ be arbitrary elements of order $> 2$. Then $ sup
(cl([ab, ba]^n)\colon n\in \mathbb{N}) = \infty.$
\end{prop}

\begin{proof}
Suppose that $\rm{sup (}cl([ab, ba]^n):\ n\in \mathbb{N}\rm{)}=k$.
 Then for all $n$ the element $[ab, ba]^n$ is a product of $\leq k$
commutators. It follows by Proposition \ref{prop:k commutators}
 that $|w_b([ab, ba]^n)| \leq 12k-3.$ On the other hand since
the word $[ab,ba]$ is cyclicaly reduced, we have  the equality
$\it{w}_b([ab, ba]^n) = -2n$ which generalizes the equality
(\ref{equation:ab}) (cf. \cite[p.579, Proof of Lemma 2.28]{Rh}  and \cite[p.578, Proof of Lemma 2.26, Case
1]{Rh}). If $n$ is large enough then $|-2n| > 12k-3$ and we get a
contradiction.
\end{proof}

\begin{rmk}
If the order of the element $b$ is 2 or, more generally,
if\ \ \ $b^2a =
ab^2$ in a group $G$ then $[a,b]^{2n} = [[a,b]^{-n}, b]$ and
$[a,b]^{2n+1} = [a[a,b]^{-n}, b],$ i.e. $cl([a,b]^{2n}) = 1 =
cl([a,b]^{2n+1}) $ for every $n.$
\end{rmk}
\begin{proof}
Since $[a,b] = b^{-1}[a,b]^{-1}b$ it follows that
$[a,b]^n = b^{-1}[a,b]^{-n}b.$
By
multiplying both sides of the equality 
by
$[a,b]^n$ we get
$[a,b]^{2n} = [[a,b]^{-n}, b].$
Since $b[a,b]^n = [a,b]^{-n}b$ and
$[a,b]^{2n+1} = [a,b]^na^{-1}b^{-1}ab[a,b]^n$
we get $[a,b]^{2n+1} = [a[a,b]^{-n},b]$.
\end{proof}

\begin{prop}\label{prop:Borsuk} {\rm (Borsuk \cite{B})}.
Consider a triple of continua $\widehat{W}\supset W \supset X$ where the
space $W$ is a strong deformation retract of $\widehat{W}$. Then
the quotient space $W/X$ is a strong deformation retract of \
$\widehat{W}/X.$ In particular, $W/X$ and $\widehat{W}/X$
then
have the same homotopy type.
\end{prop}

\section{ construction of compactum C}

Let $\mathcal{P}$  be any inverse sequence of finite polyhedra
with piecewise-linear bonding
mappings:
$$P_0 \overset {f_0} \longleftarrow P_1 \overset {f_1}
\longleftarrow P_2 \overset {f_2} \longleftarrow \cdots .$$

We denote the {\it infinite mapping
cylinder} of $\mathcal{P}$  
by $C(f_0, f_1, f_2, \dots )$ (cf. \cite{Si})
and we denote
the natural compactification
of this infinite mapping cylinder by $\underleftarrow {\lim}\
\mathcal{P}$  
with the symbol
$\widetilde{\mathcal{P}}$ 
(cf. \cite{Kr}). Let $C(f_0, f_1, f_2, \dots )^*$ be the
one-point compactification of the infinite mapping cylinder.
Denote the compactification point by $p^*.$ Obviously, the
quotient space $\widetilde{\mathcal{P}}$ by $\underleftarrow
{\lim}\ \mathcal{P}$ is homeomorphic to $C(f_0, f_1, f_2, \dots
)^*.$

\begin{prop}\label{Kras} {\rm (See  \cite{Kr})}.
If $P_0$ is a point then $\widetilde{\mathcal{P}}$ is an absolute
retract ${\rm (}AR{\rm )}$.
\end{prop}

\begin{prop}\label{prop:Equivalence} Suppose that in the inverse
sequence $\mathcal{P}$ the dimensions of all polyhedra are $\leq n$
and let  $P_0$ be a one-point space. Consider $\underleftarrow \lim\
\mathcal{P}$ as a subspace of the cube $\mathbb{I}^{2n+1}.$ Then
the quotient space $\mathbb{I}^{2n+1}$ by $\underleftarrow {\lim}\
\mathcal{P}$ is homotopy equivalent to $C(f_0, f_1, f_2, \dots
)^*.$
\end{prop}

\begin{proof}
Obviously 
$\dim \widetilde{\mathcal{P}} \leq n+1$ and
therefore $\widetilde{\mathcal{P}}$ is embeddable into
$\mathbb{I}^{2n+3}.$ According to Proposition \ref{Kras},\
$\widetilde{\mathcal{P}}$ is an $AR$ and therefore a strong
deformation retract of any $AR$ which contains it. In particular,
$\widetilde{\mathcal{P}}$ is a strong deformation retract of
$\mathbb{I}^{2n+3}.$

Applying now Borsuk Theorem ~\ref{prop:Borsuk} to the triple
$\mathbb{I}^{2n+3}\supset \widetilde{\mathcal{P}} \supset
\underleftarrow {\lim}\ \mathcal{P},$ we can
conclude that the
quotient space of $\mathbb{I}^{2n+3}$ by $\underleftarrow {\lim}\
\mathcal{P}$ is homotopy equivalent to the quotient space
$\widetilde{\mathcal{P}}$ by $\underleftarrow {\lim}\
\mathcal{P},$ i.e. to the space $C(f_0, f_1, f_2, \dots )^*.$

The homotopy type of the quotient space $\mathbb{I}^{2n+3}/
\underleftarrow {\lim}\ \mathcal{P}$ does not depend on the
embedding of $\underleftarrow {\lim}\ \mathcal{P}$ to
$\mathbb{I}^{2n+3}$, by Theorem ~\ref{prop:Borsuk} and West
and Klee
\cite{Br2}. 
Therefore, applying again  Theorem
\ref{prop:Borsuk} to the triple $\mathbb{I}^{2n+3}\supset
\mathbb{I}^{2n+1}\supset \underleftarrow {\lim}\ \mathcal{P}$,
we
can conclude that $\mathbb{I}^{2n+3}/\underleftarrow {\lim}\ \mathcal{P}$
is homotopy equivalent to $\mathbb{I}^{2n+1}/\underleftarrow
{\lim}\ \mathcal{P}$ and thus
$\mathbb{I}^{2n+1}/{\underleftarrow {\lim}\ \mathcal{P}} \simeq
C(f_0, f_1, f_2, \dots )^*.$

\end{proof}

Suppose that $P_0$ is a singleton and that  for $n > 0$ $P_n$
is a
bouquet of 4
oriented circles with the base point $p_n$. The fundamental group
$\pi_1(P_n)$ is a free group with natural generators $x_{n,1},\
x_{n,2},\ x_{n,3},\ x_{n,4}.$ Consider $\pi_1(P_n)$ as a free
product of 
free groups $F(x_{n,1};\ x_{n,2})$ and $F(x_{n,3};\
x_{n,4}).$ Let $y_{n,1},\ y_{n,2}$ and $y_{n,3}\ y_{n,4}$ be free
generators of the commutator subgroups of the groups $F(x_{n,1};\
x_{n,2})$ and $F(x_{n,3};\ x_{n,4})$, respectively. For example,
according to Proposition \ref{prop:Generators} we can suppose
that: $y_{n,1} = [x_{n,1},\ x_{n,2}], \ y_{n,2} = [x_{n,1}^2,\
x_{n,2}^2],\ y_{n,3} = [x_{n,3},\ x_{n,4}],\ y_{n,4} =
[x_{n,3}^2,\ x_{n,4}^2].$

Suppose  that $f_0$ is a trivial mapping and that
for $n
> 0$,
the mapping
$f_n:P_{n+1}\to P_n$
 is piecewise-linear and 
 such that
 $f_n(p_{n+1}) =
(p_n)$ and $f_{n\sharp}(x_{n+1,i}) = y_{n,i}$, for $i =1, 2, 3,
4,$ where $f_{n\sharp}$ is a homomorphism of the corresponding
fundamental groups induced by $f_n$. All homomorphisms
$f_{n\sharp}$ are monomorphisms since by our choice the elements
$y_{n,1},\ y_{n,2}$ and $y_{n,3},\ y_{n,4}$  are free generators.
Therefore we can consider the elements $x_{n,i}$ for $i=1,2,3,4$
as elements of the group
\begin{equation}\label{equation:Free}
F = F_{x_{1,1};\ x_{1,2};\ x_{1,3};\ x_{1,4}}.
\end{equation}

The 1-dimensional compactum $C = \underleftarrow {\lim}\
\mathcal{P}$ is embeddable into the interior $(0,1)^3$ of the cube
$\mathbb{I}^3.$ Hereafter we shall fix such an embedding.
Since $ y_{n,i}$ belongs to the commutator subgroup of $\pi_1(P_n)$
and since dim $P_n = 1$ it follows that
$f_{n\ast}:H_{\ast}(P_{n+1})\to H_{\ast}(P_{n})$ is
a
trivial
mapping. By the
Universal Coefficient Theorem the mapping
$H^{\ast}(P_{n})\to H^{\ast}(P_{n+1})$ is also trivial and the following holds:
\begin{equation}\label{equation:Cech}
\check{H}^*(C) \cong \check{H}^*(pt)
\end{equation}

\section{Proof of Main Theorem }

To prove Main Theorem 1.3 we must first prove the following:

\begin{thm}\label{thm:Homology}
The 1-dimensional singular homology group of the quotient space
$\mathbb{I}^{3}/C$ is uncountable.
\end{thm}

\begin{proof}

Let $I_n$ be the segment connecting $p_{n+1}$ and $p_n$ in the
space \\
$C(f_0,f_1,f_2,\dots)^*.$ Since $f_{n\sharp}$ maps the
groups
$F(x_{n+1,1},\ x_{n+1,2})$ and $F(x_{n+1,3},\ x_{n+1,4})$ to
$F(x_{n,1},\ x_{n,2})$ and $F(x_{n,3},\ x_{n,4}),$ respectively,
the space $C(f_0,f_1,f_2,\dots)^*$ splits, i.e. it
is the
union of two
closed subspaces, the intersection of which is the segment
$\{p^*\}\bigcup(\cup_{i=0}^{\infty}I_n).$

Suppose that $H_1(\mathbb{I}^{3}/C)$ were countable. Then
$H_1(C(f_0,f_1,f_2,\dots)^*)$ would also be countable, according
to Proposition~\ref{prop:Equivalence}. From the following
Mayer-Vietoris exact sequence:
$$H_1(P_1)\to H_1(C(f_0)) \oplus H_1(C(f_1,f_2,f_3,\dots)^*)\to
H_1(C(f_0,f_1,f_2,\dots)^*)\to 0$$
it would then follow that the group $H_1(C(f_1,f_2,f_3,\dots)^*)$ is
countable.

Define inductively a sequence of loops $g_n$ with the base point
$p_1 \in C(f_1,f_2,f_3,\dots)^*$. Let $g_1$ be any loop in $P_1$
representing a commutator of the group $\pi_1(P_1).$ Suppose that
the loops $g_i$ are defined for  $i \leq n-1$. Then define the
loop $g_n$ in the following way. Consider the set of loops
${g_1}^{\varepsilon_1}{g_2}^{\varepsilon_2}{g_3}^{\varepsilon_3}
\cdots
 {g_{n-1}}^{\varepsilon_{n-1}},$ where every $\varepsilon_i$ is
 equal to $0$ or
 $\pm 1$.
 This set is finite and therefore there exists the maximum of the
 commutator
length of its elements. Call this number $K_n.$ We have
\begin{equation}\label{equation:K_n}
cl([{g_1}^{\varepsilon_1}{g_2}^{\varepsilon_2}
{g_3}^{\varepsilon_3}\cdots {g_{n-1}}^{\varepsilon_{n-1}}]) \leq
K_n
\end{equation}

Consider the elements $x_{n, 1}, x_{n, 3}$ as elements of the
group $F$ \ (see (\ref{equation:Free})). These are non-neutral
elements of the commutator subgroup and therefore by Proposition
\ref{prop:Series}
there exists a
finite number $m_n$ \ such that $x_{n, 1}\notin F_{m_n}$ and
$x_{n, 3}\notin F_{m_n}$. Let $Y_{m_n} =
C(f_1,f_2,f_3,\dots)^*/
C(f_{{m_{n}}+1},f_{{m_{n}}+2},f_{{m_{n}}+3}\dots)^*.$ There is 
a
homomorphism $\pi_1(Y_{m_n}) \to F/F_{m_n}.$ Since by Proposition
\ref{prop:Torsions} the free nilpotent group $F/F_{m_n}$ has no
torsion, it follows that the order of any element which does not lie in the kernel
of this homomorphism 
must be infinite. Therefore the orders of the
natural elements $\widetilde{x_{n, 1}}$ and $\widetilde{x_{n,3}}$ which
correspond to $x_{n, 1}$ and $x_{n,3}$ in $\pi_1(Y_{m_n})$ are 
infinite.

The space $Y_{m_n}$ splits, therefore the group $\pi_1(Y_{m_n})$ is
a free product of two groups.
According to Proposition \ref{prop:Sup} there exists a number
$L_n$ such that
\begin{equation}\label{equation:cl>}
cl([\widetilde{x_{n,1}}\widetilde{x_{n,3}},\
\widetilde{x_{n,3}}\widetilde{x_{n,1}}]^{L_n} )>  2K_n + n.
\end{equation}

Let $g_{n}$ be a loop in $P_{m_n}\cup (\cup_{i=1}^{{m_n}-1}I_i)
\subset C(f_1,f_2,f_3,\dots)^*$ which represents the element
$[\widetilde{x_{n,1}}\widetilde{x_{n,3}},\
\widetilde{x_{n,3}}\widetilde{x_{n,1}}]^{L_n}.$ Then by an
inductive procedure we get a sequence of loops $\{g_n: n\in
\mathbb{N}\}.$

Consider now the sequence $\varepsilon$ of units and zeros:
$\varepsilon = (\varepsilon_1, \varepsilon_2, \varepsilon_3,
\cdots).$ To every such sequence there corresponds an element
$[g^{\varepsilon}] = [{g_1}^{\varepsilon_1}{g_2}^{\varepsilon_2}
{g_3}^{\varepsilon_3}\cdots]$ of the group $\pi_1
(C(f_0,f_1,f_2,\dots)^*)$ (Note that the infinite product of  loops is
not always defined but in our case obviously there exists a loop
$g^{\varepsilon}$ such that its projection to the space $Y_{m_n}$
is homotopy equivalent to the projection of
${g_1}^{\varepsilon_1}{g_2}^{\varepsilon_2}
{g_3}^{\varepsilon_3}\cdots{g_k}^{\varepsilon_k}$ for $k\ge m_n$
and the element $[g^{\varepsilon}]$ is well-defined).

Since the set of all sequences of units and zeros is uncountable
whereas
by our hypothesis the group $H_1(C(f_1,f_2,f_3,\dots)^*)$ is
countable, it follows that there exist two elements
$[g^{\varepsilon}]$ and $[g^{\varepsilon'}]$ generating the same
element in the homology group and such that ${\varepsilon_i} \neq
{\varepsilon'_i}$ for an infinite set of indices $\{i\}$. The
element $[g^{\varepsilon'}][(g^{\varepsilon})]^{-1}$ belongs to
the commutator subgroup of $\pi_1(C(f_1,f_2,f_3,\dots)^*)$ or
equivalently, it
has a finite commutator length
$cl([g^{\varepsilon'}][g^{\varepsilon}]^{-1}) = k < \infty. $ Let
$n$ be a number such that
\begin{equation}\label{equation:n>k}
n > k
\end{equation}
\noindent and $\varepsilon_n \neq \varepsilon'_n$. We can suppose
that $\varepsilon_n = 1$ and $\varepsilon'_n = 0.$ 
Consider the
mapping $\Pi_n: \pi_1(C(f_1,f_2,f_3,\dots)^*) \to \pi_1(Y_n).$ By
(\ref{equation:cl(phi)}) we have the following:
\begin{equation}\label{equation:Pr}
cl(\Pi_n([g^{\varepsilon'}][g^{\varepsilon}]^{-1})) \leq k .
\end{equation}
\noindent
 On the other hand
$$\Pi_n([g^{\varepsilon'}][g^{\varepsilon}]^{-1}) =
[{g_1}^{\varepsilon'_1}{g_2}^{\varepsilon'_2}\cdots
 {g_{n-1}}^{\varepsilon'_{n-1}}]
 {[g_n}]^{-1}[{g_{n-1}}^{-\varepsilon_{n-1}}\cdots
 {g_2}^{-\varepsilon_2}{g_1}^{-\varepsilon_1}]. $$

 By equality (\ref{equation:cl<}) we therefore have the following:
\begin{equation}
\label{equation:2cl<}
cl([g_n]) - cl([{g_1}^{\varepsilon'_1}{g_2}^{\varepsilon'_2}
 \cdots
 {g_{n-1}}^{\varepsilon'_{n-1}}]) - cl([{g_{n-1}}^{-\varepsilon_{n-1}}
 \cdots 
 {g_2}^{-\varepsilon_2}{g_1}^{-\varepsilon_1}]) 
\end{equation} 
$$ \leq  cl(\Pi_n([g^{\varepsilon'}][g^{\varepsilon}]^{-1})).
$$

 By (\ref{equation:cl>}) we have that  $cl([g_n])>  2K_n + n.$
By (\ref{equation:K_n}) we have also that
$$cl([{g_1}^{\varepsilon'_1}{g_2}^{\varepsilon'_2}\cdots
{g_{n-1}}^{\varepsilon'_{n-1}}]) \leq K_n$$ and
$$cl([{g_{n-1}}^{-\varepsilon_{n-1}}\cdots
{g_2}^{-\varepsilon_2}{g_1}^{-\varepsilon_1}])\leq K_n.$$ It
follows that the left side of the inequality (\ref{equation:2cl<})
is not less than $(2K_n + n) - K_n  - K_n = n.$ However,
by (\ref{equation:Pr})
the part on
the right does not exceed $k$. We thus get
a contradiction to (\ref{equation:n>k}).

\end{proof}

\begin{cor}\label{cor:Cohomology}
The singular cohomology group $H^*$of the space $\mathbb{I}^{3}/C$
is non-trivial.
\end{cor}

\begin{proof}
Suppose that singular cohomology groups $H^1$ and $H^2$ of the
space $\mathbb{I}^{3}/C$ were trivial. Then by the Universal
Coefficient Theorem for singular homology and cohomology it would
follow that $$\textrm{Hom}(H_1(\mathbb{I}^{3}/C), \mathbb{Z}) = 0 \ \
\hbox{\rm and} \ \
\textrm{Ext}(H_1(\mathbb{I}^{3}/C), \mathbb{Z}) = 0.$$
Therefore by the Nunke Theorem (cf. \cite[page 372, Theorem 15.6]{Br}
or \cite{Nu}) it would follow that $H_1(\mathbb{I}^{3}/C) = 0.$
This would contradict Theorem ~\ref{thm:Homology}.
\end{proof}

\begin{thm}\label{thm:Main}
The space $\mathbb{I}^{3}/C$ is a non-homologically locally
connected (non-$HLC$) 3-dimensional compact metrizable  cohomology
manifold with boundary.
\end{thm}

\begin{proof}

For $HLC$ paracompact spaces the singular cohomology and \v Cech
cohomology are naturally isomorphic \cite[p. 184, Theorem
1.1]{Br}. The \v Cech cohomology of the space $\mathbb{I}^{3}/C$
is trivial -- by the Vietoris-Begle Theorem and since the
\v Cech cohomology of the space $C$ is trivial (see
(\ref{equation:Cech})), but singular cohomology is nontrivial, as
it was shown in Corollary \ref{cor:Cohomology}. Therefore
$\mathbb{I}^{3}/C$ is not an $HLC$ space.

Other assertions of the theorem immediately follow by Wilder's
monotone mapping theorem (see \cite[p.389, Theorem 16.33]{Br}).
\end{proof}

\begin{rmk}
{\rm Theorem \ref{thm:Main} shows that the statement of the first
author in the paper \cite{K} preceding Theorem on page 531 should be
corrected as follows: 
Every compact metrizable space,
$HLC$ acyclic with respect to \v Cech
cohomology,  is an acyclic space with
respect to the singular, Borel-Moore, Steenrod-Sitnikov, Vietoris
and \v Cech homology.}

\end{rmk}

{\sl Proof of Main Theorem}: Since we can suppose that
$\mathbb{I}^n/C\subset M^n/C$ and since the boundary of
$\mathbb{I}^n/C$ is an (n-1)-dimensional sphere,  $n \geq 3,$ it
follows by the Mayer-Vietoris exact sequence that
$H_1(\mathbb{I}^n/C)$ is isomorphic to the subgroup of
$H_1(M^n/C).$

By Theorem
~\ref{thm:Homology},
the group $H_1(\mathbb{I}^n/C)$ is uncountable.
Therefore
\begin{equation}\label{equation:H}
{\text{ the group}} \ \ H_1(M^n/C)\  {\text{ is
uncountable.}}
\end{equation}

Suppose that  $H_1(M^n/C)$ were an $HLC$ space. Then its \v Cech
cohomology would be finitely generated, by Wilder's theorem (see
e.g. \cite[page 127, Theorem 17.4]{Br}). Since for $HLC$ spaces \v Cech and
singular cohomology are isomorphic, it would
follow that $H^1(M^n/C)$ and $H^2(M^n/C)$ must be finitely
generated. Then by the Universal Coefficient Formula for singular
homology and cohomology it would follow that the groups
$\textrm{Hom}(H_1(M^n/C), \mathbb{Z})$ and
$\textrm{Ext}(H_1(M^n/C), \mathbb{Z})$ are finitely generated. It would then
follow, by Bredon's Theorem
 (\cite[p. 367, Proposition 14.7]{Br}), that the group $H_1(M^n/C)$
must be finitely generated and therefore countable. This
contradicts (\ref{equation:H}). Hence the proof of Theorem 1.3
is completed. \qed

\section {Acknowledgements}
This research was supported by the
Slovenian Research Agency program
P1-0292-0101-04
and
projects J1-6128-0101 and J1-9643-0101.
We thank the referee for several useful comments and suggestions.

\end{document}